\documentclass[11pt,leqno]{article}
\usepackage{amsmath,amsthm}
\usepackage {latexsym}
\usepackage{amssymb}
\newcommand{\bea}{\begin{eqnarray}}\newcommand{\eea}{\end{eqnarray}}
\newcommand{\beq}{\begin{equation}}\newcommand{\ee}{\end{equation}}
\def\pa{\partial}

\newtheorem{theorem}{Theorem}[section]\newtheorem{lemma}[theorem]{Lemma}

\theoremstyle{remark}

\def\bm{\left( \begin{array}{cc}}\def\endm{\end{array}\right)}
\newcommand{\eq}{\end{equation}}
\def\pa{\partial}
\def \rectangle#1#2{\hbox{\vrule\vbox to #2
              {\hrule\hbox to #1{\hfil}\vfil\hrule}\vrule}}

\numberwithin{equation}{section}

\begin{document}
\author {Hans Lindblad\thanks{Part of this work was done while H.L. was a Member of the
Institute for  Advanced Study, Princeton, supported by the NSF
grant DMS-0111298 to the Institute. H.L. was also partially
supported by the NSF Grant DMS-0200226. }
\, and Avy Soffer\thanks{Also a member of the Institute of Advanced
Study, Princeton.Supported in part by NSF grant DMS-0100490.}\\
University of California at San Diego and Rutgers University}
\title {Scattering and small data completeness for the
critical nonlinear Schr\"odinger equation} \maketitle
\begin{abstract} We prove Asymptotic Completeness of one dimensional NLS with long
range nonlinearities. We also prove existence and expansion of
asymptotic solutions with large data at infinity.
\end{abstract}
\section{Introduction} We consider the problem of scattering
for the critical nonlinear Schr\"odinger equation in one space dimension:
 \beq \label{eq:nls}
 i\pa_t v+\pa_x^2 v-\beta |v|^2 v-\gamma |v|^4 v=0
 \eq
 For related results in higher dimensions see e.g.\cite{
 D1,G-V1,G-V2,G-V3,G-V4,GO}
 and cited references. In one dimension the scattering problem for
 NLS and/or Hartree long range type were studied before in
 [HKN,ST,Oz]. There are many other works in this direction, most
 are cited in the above references. The work closest to ours, as
 far as the results are concerned is \cite{HN}. In this paper the
 asymptotic completeness is proved, and the $ L^{\infty}$ decay of
 the solution is shown. We use a different method, much simpler,
 and we  get an explicit construction of the phase function and the
 asymptotic form of the solution as well. We also prove by the
 same method the existence theory of wave operators for large data
 in the repulsive case, and small data in the general nonlinear
 case. In the work of \cite{G-V4}, the Hartree equation in 3 or more
 dimensions is considered; the analysis uses, among other things
 the representation of the solution in hyperbolic coordinates which we also use. But
 the approach used in this paper is very different and much more
 involved than the work presented here. See also \cite{Nak}.

 Recall first that a solution of
the linear Schr\"odinger, i.e. $\beta=\gamma=0$, with fast decaying smooth
initial data satisfies
 \beq \label{eq:linearassumptotics}
 u(t,x)\sim t^{-1/2} e^{i x^2/4t\,}\widehat{u}(0,x/t)
 \eq
where $\widehat{u}(t,\xi)=\int u(t,x) e^{-ix \xi}\, dx$ denotes
the Fourier transform with respect to $x$ only.

 We make the following ansatz for the solution of the nonlinear problem
 \beq
 v(t,x)=t^{-1/2} e^{ix^2/4t} \, V(t,y),\qquad s=t,\quad y=x/t
 \eq
Plugging this into \eqref{eq:nls} gives, since
$(i\pa_t+\pa_x^2)\Big(t^{-1/2} e^{ix^2/4t}\Big)=0$,
 \begin{multline}
 (i\pa_t+\pa_x^2)v(t,x)
 =(i\pa_t+\pa_x^2)\Big(t^{-1/2} e^{ix^2/4t} \, V(s,y)\Big)
 =t^{-1/2}e^{ix^2/4t}\big(i\pa_t +\pa_x^2+ i (x/t)\pa_x) \big)
 V(s,y)\\
 =t^{-1/2}e^{ix^2/4t}\big(i\pa_s+s^{-2}\pa_y^2)V(s,y)
 \end{multline}
 Hence \eqref{eq:nls} becomes
 \beq\label{eq:psi}
\Psi(V)=i\pa_s V -\beta s^{-1} |V|^2 V
-\gamma s^{-2} |V|^4 V +s^{-2}\pa_y^2 V=0
 \eq

It is easy to check that the general solution to the ODE
$$
L(g)= i\,\frac{d}{d s} g-\frac{\beta}{s} |g|^2\, g-\frac{\gamma}{s^2} |g|^4\, g=0
$$
is of the form
$$
g=ae^{i\phi},\quad\text{where}\quad  \phi=-\beta
a^2\ln{|s|}\, +\gamma\frac{a^4}{s}+b
$$
for some constants $a$ and $b$.

It is therefore natural with the following ansatz for the
solution of the nonlinear problem
\beq\label{eq:defdata}
V(s,y)\sim V_0(s,y)=a(y) e^{i\phi(s,y)},\qquad
\phi(s,y)=-\beta a(y)^2 \ln{|s|}+b(y)
\eq
where $a(y)$ and $b(y)$ are any smooth sufficiently fast decaying
functions of $y=x/t$.

First we show scattering, i.e. given any $a(y)$ and $b(y)$ as
above we show that there is a solution $V$ as above.

\begin{theorem} Suppose that $a(y)$ and $b(y)$ are polynomially decaying smooth real valued functions and let $v_0(t,x)=t^{-1/2}
e^{ix^2\!/4t} V_0(t,x/t)$, where $V_0$ is given by
\eqref{eq:defdata}.
 Then if $\beta\geq 0$ or
$\beta$ is small \eqref{eq:nls} has a smooth solution $v\sim v_0$
as $t\to\infty$, satisfying
\beq
\|(v-v_0)(t,\cdot)\|_{L^\infty}+\|(v-v_0)(t,\cdot)\|_{L^2}\leq
C(1+\ln (1+t))^2 (1+t)^{-1}
\eq
\end{theorem}

We then show asymptotic completeness for small initial data, i.e.
that there is a an asymptotic expansion of the form
\eqref{eq:defdata}.

\begin{theorem} Suppose that $f\in C_0^\infty$. Then if $\varepsilon>0$ is sufficiently small \eqref{eq:nls} has a global solution with data
$v(0,x)=\varepsilon f(x)$. Moreover there are functions $a(y)$ and $b(y)$ such that with $v_0(t,x)=t^{-1/2}
e^{ix^2\!/4t} V_0(t,x/t)$, where $V_0$ is given by
\eqref{eq:defdata},  $v\sim v_0$
as $t\to\infty$;
\beq
\|(v-v_0)(t,\cdot)\|_{L^\infty}\leq
C(1+t)^{-3/2+C\varepsilon}
\eq
\end{theorem}

\section{The first order asymptotics and small data existence at infinity}
The ansatz we use is an approximate solution of the form
$$
v_0(t,x)=s^{-1/2}e^{ix^2\!/4t} V_0(s,y),\quad \text{where}\quad
V_0(s,y)=a(y) e^{i\phi(s,y)},\qquad
\phi(s,y)=-\beta a(y)^2 \ln{|s|}+b(y)
$$
where $a(y)$ and $b(y)$ are any smooth sufficiently fast decaying
functions of $y=x/t$ and $s=t$.
\begin{multline}
 \big(i\,\pa_t +\pa_x^2 -\beta |v_0|^2-\gamma|v_0|^4\big)v_0=t^{-1/2}e^{ix^2/4t}
 \Big(i\,\pa_s +s^{-2}\pa_y^2 -\beta s^{-1}|V_0|^2
-\gamma s^{-2}|V_0|^4 \Big) V_0\\
 = s^{-5/2}e^{ix^2/4t}\big( -\gamma |V_0|^4 V_0+\pa_y^2 V_0\big)=F_0.
\end{multline}
Assuming that $a(y)$ and $b(y)$ decay polynomially we have
 $$
 |\pa_s^i \pa_y^j v_0|
 \leq \frac{C_N(1+\beta\ln{|s|})^{j}}{\big(s(1+|y|)\big)^{1/2}(1+|y|)^N}
 $$
 and
\beq |\pa_s^i \pa_y^j F_0|\leq
\frac{C_N\,(1+\beta\ln{|s|})^{2+j}}{\big(s(1+|y|)\big)^{5/2}(1+|y|)^N} \eq
for any $N$. It follows that
$$
 |\pa_t^i \pa_x^j v_0|
 \leq \frac{C_N}{\big(t+|x|\big)^{1/2}(1+|x/t|)^N}
 $$
 and
\beq \label{eq:inhomdecayest}|\pa_t^i \pa_x^j F_0|\leq
\frac{C_N\,(1+\beta\ln{|t|})^{2}}{\big(t+|x|\big)^{5/2}(1+|x/t|)^N} \eq
for any $N$.

 We now
consider
$$
w=v-v_0
$$

$$
(i\,\pa_t +\pa_x^2)w= G(v_0,w)+F_0
$$
where
$$
G(v_0,w)=\beta \big( |v_0+w|^2 (v_0+w)-|v_0|^2v_0\big)
+ \gamma \big( |v_0+w|^4(v_0+w)-|v_0|^4v_0\big)
$$

The solution of the PDE \beq\label{eq:inhompde} (i\,
\pa_t+\pa_x^2)w=F \eq
with vanishing final data at infinity is
given by
$$
 w(t,x)=\int_t^\infty \int E(t-s,x-y) F(s,y)\, dy ds
 $$
 where $E$ is the forward fundamental solution of $i\pa_t+\pa_x^2$.

\begin{lemma} Suppose that
\beq\label{eq:inhoms}
i\pa_t w+\pa_x^2 w=F
\eq
Then
\beq\label{eq:inhomenergy1}
 \|w(t,\cdot)\|_{L^2}\leq \|w(t_0,\cdot)\|_{L^2}
 +\Big|\int_{t_0}^{t} \|F(s,\cdot)\|_{L^2}\, ds\Big|
 \eq
\end{lemma}
\begin{proof}
The energy identity for this equation is
\beq
\frac{d}{dt} \int |w(t,x)|^2 dx =2\int \Im{( F\overline{w})}(t,x)\, dx
\eq
where $\Im$ is the imaginary part.
\end{proof}
The energy estimate is therefore
 \beq\label{eq:inhomenergy}
 \|w(t,\cdot)\|_{L^2}
 \leq \int_t^\infty \|F(s,\cdot)\|_{L^2}\, ds
 \eq
From differentiating the equation it also follows that
 \beq\label{eq:derenergyineq}
 \sum_{|\alpha|\leq 1}\|\pa^\alpha w(t,\cdot)\|_{L^2}
 \leq \int_t^\infty \sum_{|\alpha|\leq 1}\|\pa^\alpha F(s,\cdot)\|_{L^2}\, ds
 \eq
 We will now use the above inhomogeneous estimate together with an
 iterative procedure to get existence for the equation in the
 previous section, of a solution $w$ decaying at infinity to
 zero fast, in a sense having vanishing data at infinity.
 We therefore put up an iteration
 $$
(i\pa_t+\pa_x^2)w_{0}=F_0,\qquad (i\pa_t+\pa_x^2) w_{k+1}=
G(v_0,w_k)+F_0
 $$
 where the solutions are defined as convolution with the
 fundamental solution that vanishes at infinity
 (more precise later on). We must now first find the right
 estimates for $w_0$ and thereafter make an assumption that
 the other iterates have similar bounds. It follows from
 \eqref{eq:inhomdecayest} that
 $$
 \sum_{|\alpha|\leq 1}
 \|\pa^\alpha F_0(t,x)\|_{L^2}\leq \frac{C (1+\beta\ln{|1+t|})^2}{t^{2}}
 $$
and hence
$$
\int_t^\infty \sum_{|\alpha|\leq 1}
 \|\pa^\alpha F_0(t,\cdot)\|_{L^2} dt \leq \frac{K
(1+\beta\ln{|1+t|})^2}{t}
 $$
 for some fixed constant $K$. We therefore make the inductive
 assumption that
 \beq \label{eq:indassump}
\sum_{|\alpha|\leq 1}\|\pa^\alpha w_k(t,\cdot)\|_{L^2} \leq
\frac{2K (1+\beta\ln{|1+t|})^2}{t}
 \eq
\begin{lemma}\label{decay}
\beq
\|w(t,\cdot)\|_{L^\infty}^2\leq \|w(t,\cdot)\|_{L^2} \|\pa_x w(t,\cdot)\|_{L^2}
\eq
\end{lemma}
 \begin{proof}
 Follows by H\"older's inequality $w^2\leq 2\int |w| |w_x|\, dx\leq
 2\|w\|_{L^2} \|\pa w\|_{L^2}$.
\end{proof}
It follows that
 $$
\|w_k(t,\cdot)\|_{L^\infty}\leq \frac{2K
(1+\beta\ln{|1+t|})^2}{t}
 $$
 Also using the estimates for $v_0$,
 $$
 \sum_{|\alpha|\leq 1} |\pa^\alpha v_0|\leq C_0/t^{1/2}
 $$
 we get
 $$
 \sum_{|\alpha|\leq 1} \|\pa^\alpha
 G(v_0,w_k)\|_{L^2} \leq
 \frac{C_1\beta }{t} \sum_{|\alpha|\leq 1}\|\pa^\alpha w_k\|_{L^2},
\qquad \text{if}\qquad t\geq t_0
$$
for some number $t_0=t_0(\beta)<\infty$. $t_0$ is chosen such that
the r.h.s. of equation (2.10) is smaller than $1$. Hence by the
energy inequality and the inductive assumption we get for $t\geq
t_0^\prime(\beta)$;
\begin{multline}
\sum_{|\alpha|\leq 1} \|\pa^\alpha w_{k+1}(t,\cdot)\|_{L^2}\leq
\int_t^\infty \frac{C_1\beta}{s}\sum_{|\alpha|\leq 1} \|\pa^\alpha
w_k(s,\cdot)\|_{L^2}\, ds+\frac{K(1+\ln{|1+t|})^2}{t},\\
\leq
\int_t^\infty \frac{C_1 \beta 2K (1+\beta\ln{|1+s|})^2 ds}{s^2}
+\frac{K(1+\ln{|1+t|})^2}{t}\leq \frac{(C_2\beta +1)K(1+\ln{|1+t|})^2}{t}
\end{multline}
Hence \eqref{eq:indassump} follows also for $k$ replaced by $k+1$,
if $\beta$ is so small that $C_2\beta\leq 1$.
This proves the theorem for small $\beta$.

\section{Global existence and decay for the initial value problem}
Here we show that (1.1) has a global solution for small initial
data and that the solution decays like $t^{-1/2}$. Let us suppose
we are given initial data when $s=t=1$,say in $C_0^{\infty}$ .
\begin{lemma}\label{decaylemma}  Suppose that $g$ is real valued and
\beq
i\pa_s V -g V=F.
 \eq
Then
\beq
|V(s)|\leq |V(s_0)|+\Big|\int_{s_0}^s |F(\sigma)|\, d\sigma\Big|.
\eq
\end{lemma}
\begin{proof} Multiplying with the integrating factor $e^{iG(s)}$, where
$G=\int g\, ds$ gives
$\pa_s \big( Ve^{iG}\big)=-i \, F\, e^{iG}$ and the lemma follows from
integrating this.
\end{proof}
It now follows that if \eqref{eq:psi} holds then
\beq \label{eq:inhomdecay}
|V(s,y)|\leq
|V(1,y)|+\int_1^s |\pa_y^2 V(\sigma,y)|\frac{d\sigma}{\sigma^2}.
\eq Hence the desired bound for $V$ would follow if we can prove
that for some fixed $\delta<1$; \beq\label{eq:need} |\pa_y^2
V(s,y)|\leq C\varepsilon (1+s)^{\delta}. \eq
We will now derive
this bound from energy bounds. We will assume that
\beq\label{eq:use} |V|\leq C_0 \varepsilon \eq Writing
\eqref{eq:psi} in the form
\beq\label{eq:Psi} i\pa_s V
-s^{-2}\pa_y^2 V=F=\beta s^{-1}|V|^2\,  V +\gamma s^{-2} |V|^4 V
 \eq
and differentiating the above equation with respect to $y$ gives
\beq\label{eq:Psi} \big(i\pa_s  -s^{-2}\pa_y^2 \big)
V^{(k)}=F^{(k)},\qquad V^{(k)}=\pa_y^k V,\quad F^{(k)}=\pa_y^k F
\eq We claim that \beq \|F^{(k)}(s,\cdot)\|_{L^2}\leq C
s^{-1}\big(1+s^{-1}\|V(s,\cdot)\|_{L^\infty}^2\big)
\|V(s,\cdot)\|_{L^\infty}^2 \|V^{(k)}(s,\cdot)\|_{L^2} ,\quad
k=0,1,2,3. \eq In fact
\begin{align}
|F^{(0)}|&\leq C s^{-1} (1+s^{-1}|V|^2) |V|^3\\
|F^{(1)}|& \leq C s^{-1} (1+s^{-1}|V|^2) |V|^2 |\pa_y V|\\
|F^{(2)}|&\leq C s^{-1}(1+s^{-1}|V|^2)\big(|V|^2\, |\pa_y^2 V| +|V|\, |\pa_y V|^2\big)\\
|F^{(3)}|&\leq C s^{-1} (1+s^{-1}|V|^2) \big(|V|^2|\pa_y^3 V| +|V|\, |\pa_y V|\, |\pa_y^2 V|+ |\pa_y V|^3\big)
\end{align}
For $k=0,1$ this is obvious and for $k\geq 2$ this follows from
interpolation (proved by just integrating by parts):
\begin{lemma}
\beq
\|\pa_y^j V\|_{L^{2k/j}}^{k/j}\leq C\|V\|_{L^\infty}^{k/j-1} \|\pa_y^{\,k} V\|_{L^2}
\eq
\end{lemma}
For $k=2$ we have
\beq
\|\pa_y V(t,\cdot)\|_{L^4}^2\leq
C \|V(t,\cdot)\|_{L^\infty} \|\pa_y^2 V(t,\cdot)\|_{L^2},
\eq
Similarly for $k=3$ we have
\begin{align}
\|\pa_y V(t,\cdot)\|_{L^6}^3 &\leq
C \|V(t,\cdot)\|_{L^\infty}^2 \|\pa_y^3 V(t,\cdot)\|_{L^2},\\
\|\pa_y^2 V(t,\cdot)\|_{L^3}^{3/2} &\leq
C \|V(t,\cdot)\|_{L^\infty}^2 \|\pa_y^3 V(t,\cdot)\|_{L^2}.
\end{align}
\begin{lemma}\label{Energy2} Suppose that
\beq
i\pa_s W-s^{-2}\pa_y^2 W=F.
\eq
Then
\beq\label{eq:inhomenergy1}
 \|W(s,\cdot)\|_{L^2}\leq \|W(s_0,\cdot)\|_{L^2}
 +\Big|\int_{s_0}^{s} \|F(\sigma,\cdot)\|_{L^2}\, d\sigma\Big|.
 \eq
\end{lemma}

Assuming the bound \beq \|V(s,\cdot)\|_{L^\infty}\leq
K\varepsilon=\delta\eq with a constant independent of $s$ we have
hence proven, using the above lemma and equations (3.7)-(3.16)
that
 \beq\label{eq:energyk}
\|V^{(k)}(t,\cdot)\|_{L^2}\leq \|V^{(k)}(1,\cdot)\|_{L^2} +\int_1^
t D_k(\delta+\delta^2)\|V^{(k)}(\tau,\cdot)\|_{L^2}\,\tau^{-1}  d\tau,
\eq from which it follows that \beq \|V^{(k)}(t,\cdot)\|_{L^2}\leq
C_k\varepsilon (1+t)^{D_k(\delta+\delta^2)},\qquad k=0,1,2,3. \eq and by
Lemma \ref{decay} \beq \|V^{(k)}(t,\cdot)\|_{L^\infty}\leq
 C_{k+1}\varepsilon (1+t)^{D_k(\delta+\delta^2)},\qquad k=0,1,2 \eq
where $C_k$ is a constant such that $\sum_{j\leq k}\|V^{(j)}(1,\cdot)\|_{L^2}\leq C_k\varepsilon$. Now suppose that $\varepsilon>0$ is so small that
\beq
D_3\big( 4C_3\varepsilon+(4C_3\varepsilon)^2\big)\leq 1/4.
\eq
It then follows from \eqref{eq:inhomdecay} that
\beq
\|V(s,\cdot)\|_{L^\infty}\leq 4C_3\varepsilon\qquad\text{and}\qquad
\|\pa_y^2 V(s,\cdot)\|_{L^\infty}\leq C_3\varepsilon (1+t)^{1/2}.
\eq

This is
more than needed in \eqref{eq:need}.

\section{The completeness}
\begin{lemma}\label{completelemma}  Suppose that $g$ is real valued and
\beq
i\pa_s V -g V=F
 \eq
Then with $G(s,y)=\int_{s_0}^s g(\tau,y)\, d\tau$
\beq
\big|\pa_s |V|\big| +\big|\pa_s (Ve^{iG})\big|\leq 2|F|\,
\eq
\end{lemma}
\begin{proof} Multiplying with $\overline{V}$ gives
\beq
i\pa_s |V|^2=\Im F \overline{V}
\eq
and it follows that $|\pa_s |V||\leq |F|$.
Multiplying with the integrating factor $e^{iG(s)}$, where
$G=\int g\, ds$ gives
$\pa_s \big( Ve^{iG}\big)=-i \, F\, e^{iG}$ and the lemma follows.
\end{proof}
In the application
$g=-\beta s^{-1}|V|^2
-\gamma s^{-2} |V|^4$ and $F=s^{-2}\pa_y^2 V$.
We already have proven that
\beq
|F(s,y)|\leq C\varepsilon s^{-2+C\varepsilon^2}
\eq
in the previous section. It therefore follows from the above lemma that
the limit exists
\beq
\Big| |V(s,y)|-a(y)\Big|\leq C\varepsilon s^{-1+C\varepsilon^2},
\qquad\text{where}\quad a(y)=\lim_{s\to\infty} |V(s,y)|
\eq
It therefore also follows from the lemma that
\beq
|G(s,y)-\phi(s,y)|\leq C\varepsilon s^{-1+C\varepsilon^2}, \qquad \text{where}\quad \phi(s,y)=a(y)^2 \beta \ln{|s|}+b(y)
\eq
and $-b(y)$ is defined as the limit of $G(s,y)-a(y)^2\beta\ln{|s|}$
as $s\to\infty$. Hence
\beq
\Big| V(s,y)-a(y)e^{i\phi(s,y)}\Big|\leq
C\varepsilon s^{-1+C\varepsilon^2}
\eq

\section{Higher order asymptotics and large data existence at infinity}
We now want to construct a higher order asymptotic expansion at
infinity. Therefore, we want to linearize the operator
$$
L(g)= i\,\frac{d}{d s}
g-\frac{\beta}{s}G_1(g)-\frac{\gamma}{s^2}G_2(g),\qquad
G_1(g)=|g|^2\, g,\qquad G_2(g)=|g|^4\, g
$$
We have $G_i(V_0+W)=G_i(V_0)+G_i^{\,\prime}(V_0) W+O(|W|^2)$, where
$$
G_1^{\,\prime}(V_0)W=2|V_0|^2 W+V_0^2 \overline{W},\qquad
G_2^{\,\prime}(V_0)W=3|V_0|^4 W+2|V_0|^2V_0^2 \overline{W}
$$
Here,
$$
V_0=a(y)e^{i\phi(s,y)} \qquad \phi(s,y)=-\beta a(y)^2\ln{|s|}+b(y)
$$

 Hence the linearized operator is
$$
L_0 W=L^{\,\prime}(V_0)W= i\,\frac{d}{ds} W-\frac{\beta}{s}
G_1^{\,\prime}(V_0) W-\frac{\gamma}{s^2} G_2^{\,\prime}(V_0) W
$$
Observe that $L_0$ is not complex linear. If $Z$ is constant it
therefore follows that $(k \geq 1)$
$$
L_0 \Big( \frac{e^{i\phi}\ln^{\,j}{|s|}}{s^k}Z\Big)=
\frac{e^{i\phi}\ln^{\,j}{|s|}}{s^{k+1}}\big( (-2\beta a^2 -i
k)Z-\beta a^2 \overline{Z}\big)+ i\,j
\frac{e^{i\phi}\ln^{\,j-1}{|s|}}{s^{k+1}} Z
+\frac{e^{i\phi}\ln^{\,j}{|s|}}{s^{k+2}}\big( -3\gamma a^2
Z-2\gamma a^4 \overline{Z}\big)
$$
The inverse of
$$
(-2\beta a^2 -i k)Z-\beta a^2 \overline{Z}=Y
$$
is given by
$$
Z=\frac{1}{k^2+3\beta^2a^4} \big(-2\beta a^2 +i\, k\big)
Y+\frac{1}{k^2+3\beta^2a^4} \beta a^2 \overline{Y}
$$
and hence
\begin{multline}
L_0\Big(\frac{e^{i\phi}\ln^{\,j}{|s|}}{s^k
(k^2+3\beta^2a^4)}\Big(\big(-2\beta a^2 +i\, k\big) Y+\beta a^2
\overline{Y}\Big)\Big) =\frac{e^{i\phi}\ln^{\, j}{|s|}}{s^{k+1}}Y+
i\,j\frac{e^{i\phi}\ln^{\,j-1}{|s|}}{s^{k+1}(k^2+3\beta^2a^4)}\Big(
\big(-2\beta a^2 +i\, k\big) Y
+\beta a^2 \overline{Y}\Big)\\
+\frac{e^{i\phi}\ln^{\,j}{|s|}}{s^{k+2}(k^2+3\beta^2a^4)} \gamma
a^2\big((2\beta a^2 -i\, k)[3Y+2a^2\overline{Y}]+\beta
a^2(2a^2Y-3\overline{Y})\big)
\end{multline}

It follows that
\begin{lemma} Let ${\cal S}_k$ denote a finite sum of the form $(k
\geq 1)$
 \beq
 \sum_{k^\prime\geq k,\, j\geq 0} c_{k^\prime j}(y)\frac{e^{i\phi}\, \ln^{\,
 j}{|s|}}{s^{k^{\,\prime}}},\qquad \phi=-\beta a(y)^2 \ln{|s|}+b(y),
\eq with coefficients decaying polynomially in $y$. More precisely
$|\pa^\alpha c_{j k^{\,\prime}}(y)|\leq C_N (1+|y|)^{-N}$, for any
$N$ and $c_{j k^{\,\prime}}=0$ for $k^{\,\prime}$, $j$
sufficiently large. Here $\ln^{\, 0}{|s|}=1$.

 Then if $k\geq 1$ and $\psi_{k+1}\in {\cal S}_{k+1}$ there is $\phi_k\in {\cal
 S}_k$ and $\psi_{k+2}\in {\cal S}_{k+2}$ such that
 \beq
 L_0 \, \phi_k=\psi_{k+1}+\psi_{k+2}
 \eq
\end{lemma}
Recall that $\Psi(V)=L(V)+s^{-2}\pa_y^2 V$ and that
$L_0=L^\prime(V_0)$. We have
 \begin{lemma} Let $\Psi_n =\Psi^{\,\prime}(V_n)$ and
suppose that $V_n-V_0\in {\cal S}_1$. Then if $k\geq 1$ and
$\psi_{k+1}\in {\cal S}_{k+1}$ there is $\phi_k\in {\cal
 S}_k$ and $\psi_{k+2}\in {\cal S}_{k+2}$ such that
 \beq
 \Psi_n  \, \phi_k=\psi_{k+1}+\psi_{k+2}
 \eq
\end{lemma}
\begin{proof} First, let $\phi_k\in {\cal S}_k$ be as in the previous lemma.
Then $(\Psi_0-L_0)\phi_k=s^{-2}\pa_y^2 \phi_k\in {\cal S}_{k+2}$.
Furthermore $\Psi_n -\Psi_0 =s^{-1}\beta G^{\,\prime} (V_n)
-s^{-1}\beta G^{\,\prime}(V_0)=s^{-1} O(V_n-V_0)\in {\cal S}_2$ so
$\big(\Psi_n-\Psi_0\big)\phi_k \in {\cal S}_{k+2}$.
\end{proof}

By the results of previous sections,
 $\Psi(V_0)\in {\cal S}_2$. See e.g. equation (2.1).
We will now inductively, for $n\geq 1$ construct $V_n$ such that
$V_n-V_0 \in {\cal S}_1$ and $\Psi(V_n)\in {\cal S}_{n+2}$. Assume
that this is true for $n\leq k$. Then by the above lemma ( with
$\psi_{k+1}=\Psi(V_k) \in {\cal S}_{k+2})$ we can find $V_{k+1}$
such that
 \beq
 \Psi(V_k)+\Psi^\prime(V_k)(V_{k+1}-V_k)\in {\cal S}_{k+3},
 \qquad V_{k+1}-V_k\in{\cal S}_{k+1}.
 \eq
Furthermore, there are bilinear forms in $(X,Z)$;
$G_i^{\prime\prime}(U,V)(X,Z)$ such that
 \beq
 G_i(U)=G_i(V)+G_i^{\,\prime}(V)(U-V)+G_i^{\prime\prime}(U,V)(U-V,U-V)
 \eq
Then \beq \Psi(U)=\Psi(V)+\Psi^{\,\prime}(V)(U-V)
-\frac{\beta}{s}G_1^{\prime\prime}(U,V)(U-V,U-V)
-\frac{\gamma}{s^2}G_2^{\prime\prime}(U,V)(U-V,U-V).
 \eq
 Hence
\begin{multline}
\Psi(V_{k+1})=\Psi(V_k)+\Psi^{\,\prime}(V_k)(V_{k+1}-V_k)\\
-\frac{\beta}{s}G_1^{\prime\prime}(V_{k+1},V_k)(V_{k+1}-V_k,V_{k+1}-V_k)
-\frac{\gamma}{s^2}G_2^{\prime\prime}(V_{k+1},V_k)(V_{k+1}-V_k,V_{k+1}-V_k)\in
{\cal S}_{k+3}.
 \end{multline}

Let \beq
 v_k(t,x)=t^{-1/2} e^{ix^2/4t} \, V_k(t,y),\qquad s=t,\quad y=x/t
 \eq
 Then (see equation (2.1))
\beq i\pa_t v_k+\pa_x^2 v_k-\beta |v_k|^2 v_k-\gamma |v_k|^4 v_k=
t^{-1/2} e^{ix^2/4t} \Psi(V_k)=F_k\eq
It follows that
 \beq |\pa^\alpha F_k|\leq \frac{C_k}{(t+|x|)^{2+k}}\eq
and hence
$$
 \sum_{|\alpha|\leq 1}
 \|\pa^\alpha F_N(t,\cdot)\|_{L^2}  \leq \frac{K_N
}{t^{N}}
 $$
 for some  constant $K_N$.
We then define $w_0=0$ and for $k\geq 1$: \beq
 (i\pa_t+\pa_x^2)
w_{k+1}=\beta G(v_{N},w_k)w_k+F_N,\qquad k\geq 0.
\eq

 We will
inductively assume that
 \beq\label{eq:inductiveassumptionN}
 \|\pa w_k(t,\cdot)\|_{L^2}+\|w_k(t,\cdot)\|_{L^2} \leq \frac{4K_N}{Nt^{N}}
 \eq
 Since by H\"older's inequality
 $$
 w^2\leq 2\int |w| |w_x|\, dx\leq
 2\|w\|_{L^2} \|\pa w\|_{L^2}\leq \|\pa w\|^2_{L^2}+\|w\|_{L^2}^2
 $$
  we also get
 $$
\|w_k(t,\cdot)\|_{L^\infty}\leq \frac{4K_N}{N\,t^N}
 $$
 Since also
 $$
 |V_0|= |a(y)|\leq C_0
 $$
 it follows that
 $$
\| v_N(t,\cdot)\|_{L^\infty}=t^{-1/2} \|V_N(t,\cdot)\|\leq
t^{-1/2}[\|V_0\|+C_Nt^{-1}]\leq \frac{2C_0}{t^{1/2}},\qquad t\geq
t_N=\frac{C_N}{C_0}
 $$
 since by construction, $V_N-V_0\in {\cal S}_1 $.
 So $C_0$ is independent of $N$, if $t_N$ is sufficiently
 large.
 It follows that
 \beq
 \|G(v_1,w_k)(t,\cdot)\|_{L^\infty}\leq \frac{8C_0}{t},\qquad t\geq
 t_N^{\,\prime}
 \eq
Hence by the energy inequality \eqref{eq:derenergyineq}, (5.12-14)
$$
\|\pa w_{k+1}(t,\cdot)\|_{L^2}+\|w_{k+1}(t,\cdot)\|_{L^2}\leq
\int_t^\infty \frac{\beta 8C_0}{s} \frac{ 4K_N}{N s^N}
ds+\frac{K_N}{s^{N+1}}\,
ds=\Big(\frac{32\beta\,C_0}{N}+1\Big)\frac{K_N}{N\, t^N} \leq
\frac{2K_N}{N t^N} ,\qquad t\geq t_N^{\,\prime\prime}
$$
if $\beta>0$ is sufficiently small and $t_N^{\,\prime\prime}$ is
sufficiently large. Hence \eqref{eq:inductiveassumptionN} follows
also for $k+1$.

 \end{document}